\documentclass{article}

\usepackage[english]{babel}

\usepackage[a4paper,top=2cm,bottom=2cm,left=3cm,right=3cm,marginparwidth=1.75cm]{geometry}
\usepackage[T1,T2A]{fontenc}
\usepackage[utf8]{inputenc}
\usepackage{amsmath}
\usepackage{graphicx}
\usepackage[colorlinks=true, allcolors=blue]{hyperref}

\usepackage{amsfonts} 

{\par\noindent{\bf Доказательство.}} 
{\hfill$\scriptstyle\blacksquare$} 

\begin{document}


\begin{Large}
\centerline{A note on Kronecker's approximation theorem}
\vskip+0.5cm
\centerline{D. Maksimova}
\end{Large}

\vskip+1cm

\section{Introduction}

The following statement is well known as Kronecker's approximation theorem
(see \cite{cas}, Ch III, Theorem IV).

\vskip+0.3cm
{\bf Theorem A.}
{\it  Let
\begin{equation}\label{11}
L_j(\mathbf{x})=\sum_{i=1}^m \theta_{i,j} x_i \quad  1 \leq j \leq n
\end{equation}
be $n$ homogeneous linear forms in $m$ real variables $\pmb{x}= (x_1,...,x_m)$.
Then each of the following two statements  about a real vector
$\pmb{\alpha} = (\alpha_1,...,\alpha_n)$  implies the other.

\noindent
1) For every $\varepsilon > 0 $ there exists
 $\bold{a} \in \mathbb{Z}^m $
 such that
$$
  \max_{1\le j \le n}  \| L_j(\bold{a}) - \alpha_j \| < \varepsilon,
$$

\noindent
2) If $ \bold{u} \in \mathbb{Z}^n  $ is any integer vector  such that
$$ u_1L_1(\boldsymbol{x}) + \dots + u_n L_n(\bold{x}) \in \mathbb{Z}[x_1,\dots,x_m],
$$
then
  $u_1 \alpha_1 + \dots + u_n \alpha_n \in \mathbb{Z}$.

}
\vskip+0.3cm

Everywhere in this paper $||\cdot ||$ stands for the distance to the nearest integer.

In particular in the case $m=1, N = n+1$ Theorem A  leads to the following corollary.

\vskip+0.3cm
{\bf Theorem B.}
{\it Let real numbers
$\lambda_1,...,\lambda_N$ be linearly independent    over $\mathbb{Q}$ and let
$\alpha_1,...,\alpha_N$ be arbitrary real numbers.
Then for
any $\varepsilon > 0 $ there exist arbitrary large $ t\in \mathbb{R}_+$ such that
$$
\max_{1\le j \le N} || \lambda_j t -\alpha_j||< \varepsilon.
$$}
\vskip+0.3cm

In \cite{GM} Gonek and Montgomery proved the following quantitative version of Theorem B.

\vskip+0.3cm

{\bf Theorem C.} {\it
Let $ \varepsilon_j \in \left(0,\frac{1}{2}\right), 1\le j \le N$ and
\begin{equation}\label{M}
M_j =\left\lceil \frac{1}{\varepsilon_j} \, \log \,  \frac{N}{\varepsilon_j}\right\rceil,\,\,\,\, 1\le j \le N.
\end{equation}
Suppose real numbers
$\lambda_1,...,\lambda_N$ satisfy the condition
 that
$$
\delta =
  \min_{
    \substack{(m_1,\dots,m_N) \in  \mathbb{Z}^N\setminus\{(0,...,0)\}: \\ |m_j| \leq M_j, 1\le j \le N}}
|m_1\lambda_1+...+m_N\lambda_N|>0
$$
and $\alpha_1,...,\alpha_N$ be arbitrary real numbers.
Then for  $T= \frac{4}{\delta}$ and for any $X$ there exists a real number $t$ such
that
$$
X\le t \le X+T\,\,\,\,\,\text{and}\,\,\,\,\,
\max_{1\le j \le N } ||\lambda_j t - \alpha_j || \le \varepsilon_j  .
$$
}

Gonek and Montgomery \cite{GM} used an analytic approach with involves trigonometric polynomials and goes back to Turan \cite{T} and Chen \cite{C}. Related topics were discussed recently in  \cite{KR}.
\vskip+0.3cm
In the present paper we show that the classical approach from \cite{cas} related to transference argument gives a similar result under conditions which are weaker than those from Theorem C for small $\varepsilon_j$. Namely, we prove the following

\vskip+0.3cm
{\bf Theorem 1.} {\it
Let
\begin{equation}\label{gamma}
\delta >0,\,\,\,\,
\gamma = \frac{2^{N-2}}{ N\cdot(N!)^2},  \,\,\,\,T^* = \frac{1}{\gamma\delta}
\end{equation}
and
 $\varepsilon_j \in (0, \frac{1}{2}), \,\,1\leq j \leq N$.
  Define
\begin{equation}\label{mm}
    M_j^*= \frac{1}{\gamma \varepsilon_j} \ge 1,  \,\,\,\, 1\le j\le  N.
\end{equation}
Suppose that
 real numbers
  $\lambda_1, \dots, \lambda_N$  satisfy
\begin{equation}\label{min}
   \min_{
    \substack{(m_1,\dots,m_N) \in  \mathbb{Z}^N\setminus\{(0,...,0)\}: \\ |m_j| \leq M_j^*, 1\le j \le N}} | m_1 \lambda_1 + \dots + m_{N} \lambda_{N} | \geq \delta.
\end{equation}
Then for  any collection of real numbers  $ \ \boldsymbol{\alpha} = ( \alpha_1, \dots, \alpha_{N} ) $
and for  any $ \tau \in \mathbb{R} $ there exists  real $  t \in [\tau, \tau+T^*],$  such that
\begin{equation}
  \| \lambda_jt - \alpha_j \| \leq \varepsilon_j, \text{ }1\le j \le  N.
\end{equation}
 }

\vskip+0.3cm

We see that for small  values of $\varepsilon_j$ assumption (\ref{min}) with $M^*_j$ defined in (\ref{mm}) is weaker than that with $M_j$ from (\ref{M}). Of course, constant $\gamma$ in (\ref{gamma}) is not optimal.

Some other results dealing with  classical methods from \cite{cas} related to Kronecker's theorem were discussed recently in \cite{M}.

\vskip+0.3cm

Our paper is organised as follows. In Section \ref{S2} we formulate Theorem 2 about systems of linear forms and deduce from it Theorem 1 in Section \ref{su}.
In Section \ref{S3} we formulate a general transference theorem from Cassels' book \cite{cas} and deduce from it Theorem 2.

 \section{System of linear forms} \label{S2}

 Together with the system of linear forms (\ref{11}) we consider the system of transposed forms
 $$
R_i(\mathbf{u})=\sum_{j=1}^n \theta_{i,j} u_j,\,\,\, 1\le i \le m, \,\,\,  \mathbf{u}= (u_1,...,u_n)
$$
and put $ d = m+n$.

 Theorem 1 can be deduced from a particular case of the following general statement
 which is  analogous to Theorem XVII, Ch. V from \cite{cas} and
 which we prove in Section \ref{S3}.

\vskip+0.3cm
{\bf Theorem 2. }{\it
 Let $\boldsymbol{\alpha}=\left(\alpha_1, \ldots, \alpha_n\right)\in \mathbb{R}^n$ and  $C>0, X_i>1$.

А. Inequalities
\begin{equation}\label{eq_k17_1}
\left\|L_j(\mathrm{\bold{a}})-\alpha_j\right\| \leq \varepsilon_j,
\,\,\, 1\le j \le n;\,\,\,
\left|a_i\right| \leq X_i,\,\,\,1\le i \le m
\end{equation}
have an integer solution
 $\bold{a}\in \mathbb{Z}^m$ if  the inequality
\begin{equation}\label{eq_k17_2}
\| u_1\alpha_1+...+u_n\alpha_n \| \leq \gamma_1\max \left( \max_{1 \leq i\leq m}  X_i\left\|R_i(\bold{u})\right\|,  \max_{1 \leq j\leq n} \varepsilon_j \left|u_j\right| \right)
\end{equation}
holds for all integer $\bold{u}\in \mathbb{Z}^n$  with  $\gamma_1=d$.

B. A sufficient condition for solvability of (\ref{eq_k17_1})
in integer  $\bold{a}\in \mathbb{Z}^m$ is as follows:
Inequality
(\ref{eq_k17_2}) is valid for all integer $\bold{u}\in \mathbb{Z}^n$  with  $\gamma_1=2^{d-1}(d !)^{-2}$.
}
\vskip+0.3cm
Theorem 2 gives us the following
\vskip+0.3cm
{\bf Corollary 1.}
{\it
Assume that
$$0< \varepsilon_j<\frac{1}{2}, \,\,\, M_j^{[1]} = \frac{1}{2 \varepsilon_j\gamma_1},\,\,\, 1\leq j \leq n;
\,\,\,
0 < \delta_i < \frac{1}{2} ,\,\,\,
T_i^{[1]} = \frac{2}{\gamma_1 \delta_i} ,\,\,\,  1\le  i \le m.
$$
 Assume that
\begin{equation}\label{assume}
  \min_{
    \substack{(m_1,\dots,m_n) \in  \mathbb{Z}^n\setminus\{(0,...,0)\}: \\ |m_j| \leq M_j^{[1]}, 1\le j \le n}}
 \,\,\,
    \max_{1\le i \le m}\,\,\, \delta_i^{-1}
\| R_i(\boldsymbol{m}) \| \ge 1.
\end{equation}
Then for any real vectors
 $ \boldsymbol{\alpha} = ( \alpha_1, \dots, \alpha_n )
 \in \mathbb{R}^n$ and
 $
  \boldsymbol{\tau} = (\tau_1,...\tau_m)  \in \mathbb{R}^m
  $
  there exists an integer vector $ \bold{q} = (q_1,...,q_m)$ such that
  $$
  \tau_i \le q_i\le   \tau_i + T_i^{[1]} ,\,\,\,\ 1\le i \le m ,\,\,\,\
  \text{and}\,\,\,\,\,
     \| L_j(\boldsymbol{q}) - \alpha_j \| \leq \varepsilon_j, \,\,\,\,1\le j\le n. $$}

     \vskip+0.3cm

     Proof. First of all  check (\ref{eq_k17_2}) for Part B of Theorem 2. For $\pmb{u}= \pmb{0}$ everything is clear so we assume that  $\pmb{u}\neq \pmb{0}$. If
     $$ \gamma_1  \max_{1 \leq j \leq n} \varepsilon_j \left|u_j\right| > \frac{1}{2}
     $$
     there is nothing to prove because $ ||\xi ||\le \frac{1}{2}$ for every $\xi$.
     If
         $$ \gamma_1  \max_{1 \leq j \leq n} \varepsilon_j \left|u_j\right| \le \frac{1}{2}
     ,$$
     then
     $$
     |u_j|\le \frac{1}{2\gamma_1 \varepsilon_j} = M_j^{[1]} , \,\,\,\ 1\le j \le n
     $$
     and
     by assumption (\ref{assume})  with $ X_i = \frac{T_i^{[1]}}{4}, 1\le i \le m$ we have
     $$
     \gamma_1\max_{1\le i \le m} X_i \left\|R_i(\bold{u})\right\| \ge \frac{1}{2}
     $$
     and (\ref{eq_k17_2}) follows.

     Now we
     take $\pmb{p} = (\lceil{\tau_1}+X_1\rceil ,...,\lceil{\tau_m}+X_m\rceil)\in \mathbb{Z}^m$ and apply Theorem 2 for
     $ (\alpha_1 -L_1(\pmb{p}),..., \alpha_n -L_n(\pmb{p})) $ instead of $({\alpha_1},...,\alpha_n)$. We get integer $\mathrm{\bold{a}}$ with
     $$
     |a_i| \le X_i,\,\,\, 1\le i \le m;\,\,\,\,\,
     \left\|L_j(\mathrm{\bold{a}}+\pmb{p})-\alpha_j\right\| \leq \varepsilon_j,
\,\,\, 1\le j \le n;\,\,\,
     $$
    and
    $$
    \pmb{q}
    =\mathrm{\bold{a}}
    +
\pmb{p}= (q_1,...,q_m) \in \mathbb{Z}^m
    $$
    satisfies  the desired condition because $ \delta_i \le \frac{1}{2}$ and so
    $\frac{T_i^{[1]}}{2} \ge 1$ .
     $\Box$

\section{Corollaries for simultaneous approximations}\label{su}

Being reformulated in the case of simultaneous approximation to
real numbers ($m=1$) Corollary 1 looks as follows.

\vskip+0.3cm
{\bf Corollary 2.}
{\it
Let $\theta_1,...,\theta_n$ be real numbers.
Assume that
$$0< \varepsilon_j <\frac{1}{2}, \,\,\, M_j^{[1]} = \frac{1}{2 \varepsilon_j\gamma_1},\,\,\, 1\leq j \leq n;
\,\,\,
0 < \delta < \frac{1}{2} ,\,\,\,
T^{[1]}  = \frac{2}{\gamma_1 \delta}.
$$
 Assume that
\begin{equation}\label{assume1}
\min_{
    \substack{(m_1,\dots,m_n) \in \mathbb{Z}^n\setminus \{(0,...,0)\}:
     \\ |m_i| \leq M_j^{[1]} , 1\le j \le n}}\,\,\,
\| m_1\theta_1+...+m_n \theta_n \| \ge\delta.
\end{equation}
Then for any real vectors
 $ \boldsymbol{\alpha} = ( \alpha_1, \dots, \alpha_n )
 \in \mathbb{R}^n$ and
 $
 {\tau} \in \mathbb{R}
  $
  there exists $ {q} \in \mathbb{Z}$ such that
  $$
  \tau \le q \le \tau+ T^{[1]}, \,\,\,\,\
  \text{and}\,\,\,\,\,
     \| \theta_j q - \alpha_j \| \leq \varepsilon_j, \,\,\,\,1\le j \le n. $$}

     \vskip+0.3cm

     Now we deduce Theorem 1 from Corollary 2.

\vskip+0.3cm

     Proof. Let $\theta_i =\frac{\lambda_i}{\lambda_{N}}$, $1 \leq i \leq n$. Without loss of generality we may assume that  $\frac{|\lambda_N|}{\varepsilon_N} \geq \frac{|\lambda_i|}{\varepsilon_i}$, $ 1 \leq i \leq n$.

From the assumptions of Theorem 1 we see that
\begin{equation}\label{cor_th17_gon_proof_1}
    \min_{
    \substack{(m_1,\dots,m_n) \neq (0,\dots,0) \\ |m_i| \leq M_i^{[1]}, 1\le i \le n}}\| m_1 \theta_1 + \dots + m_n \theta_n \| > \delta_0 = \min \left( \frac{\delta}{|\lambda_{N}|}, \frac{1}{2} \right).
\end{equation}
Indeed, since $\| m_1 \theta_1 + \dots + m_n \theta_n \| = |m_1 \theta_1 + \dots + m_n \theta_n + m_N|$, where $-m_N$ is the nearest integer to $m_1 \theta_1 + \dots + m_n \theta_n$, we have
 $$|m_N|
 \leq \frac{1}{2} + \sum_i |\theta_i| |M_i|
  \leq \frac{1}{2} + \sum_i \frac{|\lambda_i|}{|\lambda_N|} \frac{1}{2 \gamma_1 \varepsilon_i}
  \leq \frac{1}{2} + \frac{n}{2 \gamma_1 \varepsilon_N} \leq \frac{n+1}{2\gamma_1 \varepsilon_N} = \frac{1}{4 \gamma \varepsilon_N} < M_N^*.
$$
Here we use that $\gamma_1 = \frac{2^{N-1}}{N!^2} \leq 1$, $M_j^{[1]}=\frac{1}{4 \gamma_1 \varepsilon_j}=\frac{1}{2N} M^*_j < M^*_j$, $1 \leq j \leq n$.

Now (\ref{cor_th17_gon_proof_1})  follows from (\ref{min}) and we can apply  Corollary 2
  with $\gamma_1 = 2\gamma N$, and $\delta_0 = \frac{\delta}{\lambda_N}$. So we   get the following statement.
For any real vector
$ \boldsymbol{\beta}$ and for any
$
\tau' \in \mathbb{R}$
there exists $q \in  [\tau', \tau'+T^{[1]}]\cap \mathbb{Z}$,   such that
$$
      \left|\left| \frac{\lambda_i}{\lambda_N} q - \beta_i \right|\right| \leq \varepsilon_i, \text{ } 1 \leq i \leq n
\,\,\,\text{(here}\,\,\, T^{[1]} = \frac{2}{\gamma_1 \delta_0}=\frac{2\lambda_N}{2\gamma N \delta}=\frac{\lambda_N}{N} T^* \text{)}.
$$
Now let $\tau = \frac{\alpha_N+\tau'}{\lambda_N}$, $t = \frac{\alpha_N+q}{\lambda_N} \in \mathbb{R}$, $t \in [ \frac{\alpha_N+\tau'}{\lambda_N},  \frac{\alpha_N+\tau'+T^{[1]}}{\lambda_N} ] \subset [\tau, \tau+T^*]$, and $\alpha_i = \beta_i + \frac{\lambda_i \alpha_N}{ \lambda_N}$. Then
$$
    \| \lambda_i t - \alpha_i \| =      \left|\left| \frac{\lambda_i}{\lambda_{N}} q - \beta_i \right|\right| \leq \varepsilon_i, \text{ } 1 \leq i \leq n
;\,\,\,\,\,
    \left|\left| \lambda_{N} t - \alpha_{N} \right|\right| = ||q|| = 0 < \varepsilon_{N}
$$
and everything is proven.
$\Box$.

\section{Application of general result}  \label{S3}

Now we formulate a general transference result which  is given in \cite{cas} as Theorem XVI, Ch. V.

\vskip+0.3cm

{\bf Theorem C.} {\it
Let
 $f_k(\bold{z})$, $g_k(\bold{w})$, $1 \leq k \leq d$
 be linear forms in variables
  $\bold{z} = (z_1, \dots , z_d)$, $\bold{w} = (w_1, \dots , w_d)$ respectively. Suppose that

$$
\sum_{k=1}^d f_k(\bold{z}) g_k(\bold{w}) = \sum_k z_k w_k
$$
holds identically.
Let  $\bold{\beta} = (\beta_1, \dots, \beta_d)\in \mathbb{Z}^d$.

А.  A necessary condition  that
\begin{equation}\label{beta}
\left|\beta_k-f_k(\mathbf{b})\right| \leq 1,  \quad(1 \leq k \leq d) ,
\end{equation}
for some integer
$\bold{b}\in \mathbb{Z}^d$ is that
$$
\left\|\sum_{k=1}^d g_k(\mathbf{w}) \beta_k\right\| \leq d \max_{1\le k \le d} \left|g_k(\mathbf{w})\right|
$$
for all   $\bold{w}\in \mathbb{Z}^d$.

В.
A sufficient condition that
 (\ref{beta}) holds for some integer $\bold{b}\in \mathbb{Z}^d$, is that
$$
\left\|\sum_{k=1}^d  g_k(\mathbf{w}) \beta_k\right\| \leq 2^{d-1}(d!)^{-2} \max \left|g_k(\mathbf{w})\right|
$$
holds for all  $\bold{w}\in \mathbb{Z}^d$.
}

\vskip+0.3cm

Proof of Theorem 2. Analogously to \cite{cas}, Theorem 2 is a particular case of Theorem C for the forms
 $$
   f_k(\mathbf{z})= \begin{cases}\varepsilon_k^{-1}\left(L_k(\mathbf{x})+y_k\right) & \text { for } k \leq n \\
X_{k-n}^{-1} x_{k-n} & \text { for } n<k \leq d\end{cases} ,\,\,\,\,\,
  g_k(\mathbf{w})= \begin{cases} \varepsilon_k u_k & \text { for } k \leq n \\
X_{k-n}\left(v_{k-n}-R_{k-n}(\mathbf{u})\right) & \text { for } n<k \leq d\end{cases}
$$
in variables
 $$
  \mathbf{z}=(\mathbf{x}, \mathbf{y})=\left(x_1, \ldots, x_m, y_1, \ldots, y_n\right), \,\,\,
 \mathbf{w}=(\mathbf{v}, \mathbf{u})=\left(v_1, \ldots, v_m, u_1, \ldots, u_n\right)
 $$
 and
 $$\pmb{\beta} = (\beta_1,...,\beta_d)  =( {\varepsilon_1}^{-1} \alpha_1, \dots , \varepsilon_n^{-1} \alpha_n,\underbrace{0,...,0}_{m}).\hskip+4cm \Box
 $$

\end{document}